\documentclass[a4paper,11pt]{article}
\usepackage{paralist}   

\usepackage{amsmath}
\usepackage{amssymb}
\usepackage{amsfonts}
\usepackage{color}

\usepackage{verbatim}

\usepackage{graphicx}
\usepackage{subfig}
\usepackage{bm}
\usepackage{tikz}

\usepackage[bookmarks,bookmarksopen=false,% Klappt die Bookmarks in Acrobat aus
                  pdfauthor={Autor},%
                  pdftitle={Titel},%
                  colorlinks=true,%
                  linkcolor=blue,%
                  citecolor=blue,%
                  urlcolor=blue,%
                ]{hyperref}

\newtheorem{theorem}{Theorem}
\newtheorem{definition}{Definition}
\newtheorem{proposition}[theorem]{Proposition}
\newtheorem{corollary}[theorem]{Corollary}

\newtheorem{lemma}[theorem]{Lemma}
\newtheorem{remark}[theorem]{Remark}
\newtheorem{example}{Example}
\DeclareMathOperator{\esssup}{ess\,sup}

\DeclareMathOperator{\modi}{ mod\ 1}
\newcommand{\proof}{\textbf{Proof:\ }}
\newcommand{\pbox}{\hfill$\Box$\\}

\newcommand{\R}{\mathbb{R}}
\newcommand{\N}{\mathbb{N}}
\newcommand{\C}{\mathbb{C}}
\newcommand{\Z}{\mathbb{Z}}

\newcommand{\I}{\mathcal{I}}

\newcommand{\F}{\mathcal{F}}
\newcommand{\G}{\mathcal{G}}

\newcommand{\To}{\mathbb{T}}
\renewcommand{\H}{\mathcal{H}}
\renewcommand{\S}{\mathcal{S}}

\newcommand{\ran}{{\sf Ran}\,}
\newcommand{\Ker}{{\sf Ker}\,}

\newcommand{\norm}[2]{\left\| #2 \right\|_{#1}}
\newcommand{\ip}[2]{ \langle {#1} ,{#2}  \rangle}
\def\V{{\mathcal V}}

\begin{document}

\title{Reproducing pairs and Gabor systems at critical density}
\date{}
\author{Michael Speckbacher and Peter Balazs\footnote{Acoustics Research Institute, Austrian Academy of Sciences,  Wohllebengasse 12-14,
1040 Vienna, Austria,
E-mail: speckbacher@kfs.oeaw.ac.at, peter.balazs@oeaw.ac.at}}
\maketitle

\begin{abstract}
 We use the concept of reproducing pairs to study Gabor systems at critical density. First, we present a generalization of the Balian-Low theorem to the reproducing pairs setting. Then, we prove our main result that there exists a reproducing partner for the Gabor system of integer time-frequency shifts of the Gaussian. In other words, the coefficients for this Gabor expansion of a square integrable function  can be calculated using inner products with an unstructured family of vectors in $ L^2(\R)$. This solves the possibly last open question for this system.
\end{abstract}

{\bf MSC2010:} 42C15, 42C40\\
{\bf Keywords:} Gabor systems, reproducing pairs, critical density, Zak transform, Balian-Low theorem  \\
%{\bf Submitted to:} Transactions of the AMS

\section{Introduction}

\begin{comment}
Frames have been introduced  in 1952 \cite{duscha52} as a generalization of orthonormal bases in the context of nonharmonic analysis. More precisely, a frame is a family of vectors in a Hilbert space that allows for stable expansions and coefficient calculation. Hence, frame theory is an important tool in many different fields from signal processing to quantum mechanics.  
However, there exist many examples of complete vector families that do not generate a frame. In order to work with such systems, different approaches have been proposed. One of which is the concept reproducing pairs \cite{spexxl14}, where one uses two families (not necessarily frames) instead of a single one that construct a bounded and invertible analysis/synthesis process.
\end{comment}

The main objective of Gabor analysis is to understand the conditions and obstructions on the family $G(g,\Lambda):=\{T_{\lambda_1}M_{\lambda_2}g\}_{\lambda\in\Lambda}\subset L^2(\R)$ to be a frame.
There exists, however, a great abundance of windows $g$ and lattices $\Lambda$ generating Gabor families which are, on the one hand, complete and, on the other hand, violate at least one of the frame bounds. The well-known Balian-Low theorem, for example, states that the window function of a Gabor frame at critical density ($\Lambda=a\Z\times a^{-1}\Z$) cannot be well localized on the time-frequency plane. In fact, there are many more properties of a window function that prevent a system from being a Gabor frame, see \cite{groe14} for an overview. It is therefore reasonable to change perspective and apply approaches beyond frame theory to Gabor families at critical density.

Several generalizations of frames, such as semi-frames or reproducing pairs have been introduced. A reproducing pair \cite{spexxl14} consists of two vector families (not necessarily frames) instead of a single one that generates a bounded and invertible analysis/synthesis process in a Hilbert space. Observe that, for Gabor families, there is a conceptual similarity to weakly dual Gabor systems \cite{fezi98}, where the analysis/synthesis process is considered in terms of a Gelfand triplet. 

The main incentive for this paper is to investigate whether the obstructions for Gabor frames at critical density still hold for reproducing pairs. That is, is there a reproducing pair where one of the two families is a Gabor system generated by a well localized window?
First, we will consider the case of reproducing pairs consisting of two Gabor systems and derive a Balian-Low like result.
Then we will turn our focus on the study of the Gabor family of integer time-frequency shifts of the Gaussian $\varphi$ which is probably the most studied object in Gabor analysis.

 Already in 1932,  von Neumann \cite{neu55} claimed without proof that the system  $\G:=\{T_k M_l\varphi\}_{k,l\in\Z}$ is complete in  $L^2(\R)$. It was only in the 1970's that
Perelomov \cite{pe71},
Bargmann et al. \cite{babugikl71} and
Bacry et al. \cite{bagrza75} presented rigorous proofs of the completeness in the context of coherent states. 

A second problem was formulated by Gabor \cite{ga46} in 1946 when he asked if 
there exists a linear coefficient map $A: L^2(\R)\rightarrow \C^{\Z\times\Z}$, such that the expansion
\begin{equation}\label{quest-eq}
f=\sum_{k,l\in\Z}(Af)(k,l)T_k M_l\varphi,
\end{equation}
converges for all $f\in L^2(\R)$. The answer to this question is more subtle. Janssen \cite{ja81} showed that such a coefficient map exists. The coefficients, however, grow polynomially and \eqref{quest-eq} converges only in the sense of tempered distributions.

As the Balian-Low theorem tells us that $\G$ cannot be a frame one might ask if there is a dual window $\gamma\in L^2(\R)$, such that $(Af)(k,l)=
\ip{f}{T_k M_l\gamma}$ in \eqref{quest-eq}.
The answer is "no". This can easily be seen using Zak transform methods which results in Bastiaans' dual window \cite{ba80}, a bounded function that is not in $L^p(\R)$ for any $1\leq p<\infty$. Moreover, it also follows from the Balian-Low like result that we will show in Section 3 or from a result by Daubechies and Janssen \cite{dauja93}.

In this paper, we will consider a problem which is intermediate to the second and third question and that is yet unsolved: can the coefficient map $A$ be calculated using inner products with an unstructured family in $L^2(\R)$, that is,
is there a system $\Psi=\{\psi_{k,l}\}_{k,l\in\Z}$, such that
$(Af)(k,l)=\ip{f}{\psi_{k,l}}$ 
with the series \eqref{quest-eq} converging weakly?

We will use a characterization of  reproducing pairs from \cite{ansptr15} in our proof to show our main result in Theorem \ref{main-repr}: the existence of a reproducing partner. This family of vectors however is totally unstructured and cannot be represented as a shift-invariant system.
Our result can be  reformulated in several contexts. For example, there is a dual system for the complete Bessel-sequence $\G$, or there is a family of vectors $\Psi$ making $(\G,\Psi)$ a reproducing pair. \\

This paper is organized as follows. In Section \ref{sec:prel} we present the basics of reproducing pairs and Gabor theory needed in the course of this article. Section \ref{sec:zak-trafo} is devoted reproducing pairs of two Gabor systems and a generalization of the Balian-Low theorem.
In Section \ref{sec:quest} we investigate the existence of a reproducing partner for the system $\G$.

\section{Preliminaries and notation}\label{sec:prel}
Let  $\omega\in\big[-\frac{1}{2},\frac{1}{2}\big)^d$ and $k\in\Z^d$.
We denote the Fourier transform of a function $f:\big[-\frac{1}{2},\frac{1}{2}\big)^d\rightarrow \C$ by
$$
\F(f)[k]=\int_{\big[-\frac{1}{2},\frac{1}{2}\big)^d}f(\omega)e^{-2\pi i k\cdot\omega}d\omega,
$$
and the discrete-time Fourier transform of a sequence $c:\Z^d\rightarrow \C$ by
$$
\F_d(c)(\omega)=\sum_{k\in\Z^d}c[k]e^{-2\pi i k\cdot\omega}.
$$
\subsection{Frames and semi-frames}\label{sec:frames}
Frames have been introduced in  1952 by Duffin and Schaeffer 
\cite{duscha52} as a generalization of orthonormal bases. A countable family of vectors $\{\psi_k\}_{k\in\I}$ in a separable 
Hilbert space $\H$ is 
called a  frame if there exist constants $A,B>0$ such that
\begin{equation}\label{framedef}
A\norm{}{f}^2 \leq\sum_{k\in\I}|\ip{f}{\psi_k}|^2\leq B\norm{}{f}^2,\ \forall\ f\in\H.
\end{equation}
For a thorough introduction to frame theory, see \cite{christ1}. Frame theory has proven its usefulness in 
many different fields like sampling theory \cite{algro01} or theoretical physics \cite{alanga00}.
However, it is sometimes impossible to satisfy both frame bounds at the same time. We will see an example of this 
situation in Section \ref{sec:gabor}. 
Hence, several generalizations of frames like reproducing pairs (see
Section \ref{sec:reppair}) or semi-frames \cite{jpaxxl09,jpaxxl12} have been introduced.
 The basic idea of semi-frames is to consider complete families $\{\psi_k\}_{k\in\I}$ that only satisfy  one of the 
 inequalities in \eqref{framedef}. In particular, $\{\psi_k\}_{k\in\I}$  is called a 
lower semi-frame if
\begin{equation}
A\norm{}{f}^2 \leq\sum_{k\in\I}|\ip{f}{\psi_k}|^2,\ \forall\ f\in\H, 
\end{equation}
and is called an upper semi-frame if
\begin{equation}
0<\sum_{k\in\I}|\ip{f}{\psi_k}|^2\leq B\norm{}{f}^2,\ \forall\ f\in\H.
\end{equation}
An upper semi-frame is often also called a complete Bessel sequence.
Many results from frame theory can be extended to the setting of semi-frames.

The following Lemma can be found in \cite[Lemma 2.5]{jpaxxl09}.
\begin{lemma}\label{semi-frames-duality}
 Let $\Phi=\{\phi_k\}_{k\in\I}$ be an upper semi-frame with bound B, and $\Psi=\{\psi_k\}_{k\in\I}$
 be a family of vectors satisfying
 $$
 \ip{f}{g}=\sum_{k\in\I}\ip{f}{\phi_k}\ip{\psi_k}{g},\ \ \forall\ f,g\in\H,
 $$
 then $\Psi$ is a lower semi-frame with lower bound $B^{-1}$.
\end{lemma}

\subsection{Reproducing pairs}\label{sec:reppair}

The concept of reproducing pairs has been introduced recently in \cite{spexxl14} and  studied in more detail in
\cite{ansptr15}. The main idea is to omit both frame bounds and to consider two vector families (instead of a single one)
that  generate a bounded and boundedly invertible analysis/synthesis process.

Although the general definition in \cite{ansptr15} is given with respect to arbitrary Borel measures, 
we will only present the discrete setting here as we will study exclusively discrete Gabor systems throughout this paper.
\begin{definition}\label{rep-pair-definition}
 Let  $\Psi=\{\psi_k\}_{k\in\I},\ \Phi=\{\phi_k\}_{k\in\I}$ be two families in $\H$.
 The pair $(\Psi,\Phi)$ is called a reproducing pair for $\mathcal{H}$ if the
 operator $S_{\Psi,\Phi}:\mathcal{H}\rightarrow \mathcal{H}$, weakly given by
 \begin{equation}\label{rep-pair-def-dis}
 \ip{S_{\Psi,\Phi}f}{g}:=\sum_{k\in\I}\ip{f}{\psi_k}\ip{\phi_k}{g},
 \end{equation}
is bounded and boundedly invertible.\\
Given a family  $\Psi$, any system $\Phi$ for which $(\Psi,\Phi)$ is a reproducing pair is called a
 reproducing partner for $\Psi$.
\end{definition}

Let $\V_\Phi(\I)$ be the space of all  sequences  $\xi  : \I \to \C$ such that
$$
\Big|  \sum_{k\in\I}  \xi[k]  \ip{\phi_k}{g} \Big| \leq c \norm{}{g},  \forall\, g \in \H.
$$
Hence, Riesz representation theorem guarantees that the 
synthesis operator $D_\Phi:\V_\Phi(\I)\rightarrow\H$, weakly given by
$$
\ip{D_\Phi \xi}{g}= \sum_{k\in\I}  \xi[k]  \ip{\phi_k}{g},
$$
is well-defined. By definition, $\V_\Phi(\I)$ is the most general domain such that the synthesis operator converges weakly.
The proof of the following result can be found in \cite[Theorem 4.1]{ansptr15}. It answers the question of 
the existence of a reproducing partner for a given family in a Hilbert space. 
 \begin{theorem}\label{theo-partner}
Let $\Phi=\{\phi_k\}_{k\in\I}\subset\H$ be a family of vectors and $\{e_k\}_{k\in\I}$ be an orthonormal 
basis for $\H$. There exists another family $\Psi$, such that $(\Psi,\Phi)$ is a reproducing pair if, and only if,
\begin{enumerate}[(i)]
\item $\ran  D_\phi =\H$ 
and \item there exists a family $\{\xi_k\}_{k\in\I}\subset \V_\Phi(\I)$ such that 
\begin{equation}\label{second-assumption}
D_\Phi \xi_k= e_k,\ \forall\ k\in\I,\hspace{0.5cm} \text{and} \hspace{0.5cm}
\sum_{k\in\I}|\xi_k[n]|^2<\infty,\ \forall\ n\in\I.
\end{equation}
\end{enumerate}
A possible reproducing partner $\Psi=\{\psi_k\}_{k\in\I}$ is then given by 
 $$\psi_n:=\sum_{k\in\I}\overline{\xi_k[n]}e_k.$$
 \end{theorem}
 
The conditions $(i)$ and $(ii)$ can be interpreted in several ways. First, Property $(i)$ ensures the existence of a linear 
operator $A:\H\rightarrow \V_\Phi(\I)$ satisfying $f=D_\Phi A(f)$, for every $f\in\H$. For an example of a complete system that does
not satisfy $(i)$, see \cite[Section 6.2.3]{ansptr15}.
Property $(ii)$ guarantees that  $A(f)$ can be calculated by taking inner products of $f$ with another family
$\Psi\subset \H$.	

Second, $(i)$ and $(ii)$ guarantee that $\{\xi_k\}_{k\in\I}$ is an orthonormal basis of its closed linear 
span with respect to the inner product $\ip{\xi}{\eta}_\Phi:=\ip{D_\Phi\xi}{D_\Phi\eta}$. The second condition of 
$(ii)$ then assures that this space is a reproducing kernel Hilbert space.

\subsection{Gabor analysis}\label{sec:gabor}
Let $\lambda=(x,\omega)\in\R^2$ be a point on the time-frequency plane, a time-frequency shift of a function $g$ by $\lambda$ is 
given by
$$
\pi(\lambda)g(t):=T_xM_\omega g(t) =e^{2\pi i \omega (t-x)}g(t-x).
$$
The short-time Fourier transform $V_g$ of a function $f$ is  given by $$
V_gf(\lambda):=\ip{f}{\pi(\lambda)g}.
$$
A Gabor system is a discrete family of functions generated by time-frequency shifts  of a single window function
$g\in L^2(\R)$ 
$$
G(g,\Lambda):=\{\pi(\lambda)g\}_{\lambda\in\Lambda}.
$$
 Let $a,b>0$, the Gabor system using a  rectangular lattice is given by
$$G(g,a,b):=\{\pi(an,bm)g\}_{n,m\in\Z}.$$
The product $(ab)^{-1}$ is called the density or redundancy of the Gabor system. If $G(g,a,b)$ is a frame, then $ab\leq1$ necessarily holds, see \cite[Corollary 7.5.1]{groe1}. 
The case $ab=1$ is called critical density. To keep notation simple, we will write for short
$
g_{an,bm}:=T_{an}M_{bm}g.
$\\

The Balian-Low theorem states that, at critical density, there are no Gabor frames using a window which is well-localized
both in time and frequency. 
\begin{theorem}[Amalgam Balian-Low Theorem]
Let $ab=1$ and the Gabor system $G(g,a,b)$ be a frame, then both $g\notin W_0(\R)$ and 
 $\hat g\notin W_0(\R)$, where $$W_0(\R):=\{f\in C(\R):\ \sum_{n\in\Z}\esssup_{x\in [0,1]}|f(x+n)|<\infty\}.$$
\end{theorem}
This obstruction motivates our approach to use reproducing pair methods to study Gabor systems generated by well localized windows at critical density.  First, we will investigate if it is possible to choose another Gabor system as the reproducing partner. Such a family  necessarily satisfies the lower but not the upper frame bound by Lemma \ref{semi-frames-duality}. Thus, the generating window function cannot be well localized. Second, we will use Theorem \ref{theo-partner} to check the existence of unstructured reproducing partners.

The analysis of the conditions of Theorem \ref{theo-partner} heavily depends  on the particular window function. Thus, we will  focus on the system of integer time-frequency shift ($a=b=1$) of the
normalized Gaussian
$$\varphi_\sigma(t)=(2/\sigma)^{1/4}e^{-\pi t^2/\sigma}.$$
However, we are convinced that the recipe for our proof also works for other window functions.
As in \cite[Section 2.2]{ja81a} we will assume that $\sigma=1$
and use the notation $\varphi:=\varphi_1$ and
$
\G:=G(\varphi,1,1).
$\\

We end this introductory part to Gabor system with defining  the modulation spaces $M_s^p$. Let  $v_s(x,\omega):=(1+|x|+|\omega|)^s$, $s>0$, and $g$ be some nonzero Schwartz function, then $M_s^p$ is defined as follows
$$
M_s^p:=\big\{f\in L^2(\R):\ V_gf\cdot v_s\in L^p(\R^2)\big\}.
$$
In particular, the spaces $M^1_s$  are commonly seen as the appropriate class of window functions for Gabor analysis. For an overview on modulation spaces, see \cite{hgfei06}.

\section{A Balian-Low like theorem for reproducing pairs}\label{sec:zak-trafo}
The Zak transform of a function $f$ is given by
 $$Z_a f(x,\omega):=\sum_{k\in\Z}f(x-a k)e^{2\pi ia\omega k}.$$
 We list some of its most important properties. An extensive description can be found for example in \cite[Chapter 8]{groe1}.
 It is easy to see that $Z_af$ is periodic in the second variable and 
 quasiperiodic in the first variable, that is, for $k,l\in\Z$ one has
 $$
 Z_af(x+a k,\omega+l/a)=e^{2\pi ia k\omega}Zf(x,\omega).
 $$
Moreover, $Z_a:L^2(\R)\rightarrow L^2(\mathcal{Q}_a)$ is a unitary operator, where $\mathcal{Q}_a:=[0,a]\times[0,1/a]$. Regularity of a function is preserved by the Zak transform. To be more precise, if
 $f\in W_0(\R)$, then $Z_af\in C(\R^2)$ and $f\in \S(\R)$ if, and only if, $Z_af\in C^\infty(\R^2)$. The most important property of the Zak transform in the context of Gabor analysis is, however, that it diagonalizes the Gabor frame operator. For  $ab=1$, $f,g,\gamma\in L^2(\R)$, it holds
\begin{equation}\label{zak-diagonal}
Z_a(S_{g,\gamma} f)=a\cdot\overline{Z_a g}\cdot Z_a\gamma\cdot Z_af.
\end{equation}
Hence, $(G(g,a,1/a),G(\gamma,a,1/a))$ is a reproducing pair if, and only if, 
\begin{equation}\label{zak-cond-rep-pair}
0<m\leq|Z_a g\cdot Z_a\gamma|\leq M<\infty, \mbox{ almost everywhere},
\end{equation}
       and $S_{g,\gamma}=I$ if, and only if, $\overline{Z_a g}\cdot Z_a\gamma=1/a$,  almost everywhere.
       
There is a close connection to Gabor Schauder bases. Heil and Powell \cite[Theorem 5.10]{hepo06} have shown that $G(g,1,1)$ is a Gabor Schauder basis if, and only if, there exists $C>0$ such that for any intervals $I,J\subset \R$ one has
$$\frac{1}{|I|^2|J|^2}\cdot\int_{I\times J}|Zg(x,\omega)|^2dxd\omega\cdot\int_{I\times J}\frac{1}{|Zg(x,\omega)|^2}dxd\omega\leq C.$$
In particular, any Schauder basis  $G(g,1,1)$ and its dual basis $G(\gamma,1,1)$ form a reproducing pair.
       
       \begin{comment}
       The proof of the Amalgam Balian-Low theorem is based on the observation that, if $Z_ag$ is continuous on $\R^2$, then $Z_ag$ has a zero in $\mathcal{Q}_a$. 
 \end{comment}            
 \begin{example}\label{ex} We try to "trick" the Amalgam Balian-Low theorem by constructing a reproducing pair using window functions $g,\gamma$ such that $g\in W_0(\R)$, that is, $g$ is well localized on the time-frequency plane. Observe that $\gamma$ is then necessarily badly localized.
 Define $\vartheta(t):=t^{1/4}(1-t)^{1/4},\ t\in[0,1],$ and 
$$
Z_ag_a(x,\omega):=e^{2\pi i x/a\cdot (a\cdot\omega \modi)}\vartheta(a\cdot\omega \modi),
$$
then $Z_ag_a$ is quasiperiodic in $x$, periodic in $\omega$ and continuous on $\R^2$. Moreover, $g_a\in W_0(\R)$, $\widehat g_a(\omega)=\vartheta(a\omega) \cdot\chi_{[0,1/a]}(\omega)\in W_0(\R)$ and $\gamma_a:=Z_a^{-1}(1/\overline{Z_ag})\in L^2(\R)$.
Hence, $\big(G(g_a,a,1/a),G(\gamma_a,a,1/a)\big)$ is a reproducing pair for every $a>0$.
\begin{figure}[t]
\includegraphics[scale=0.16]{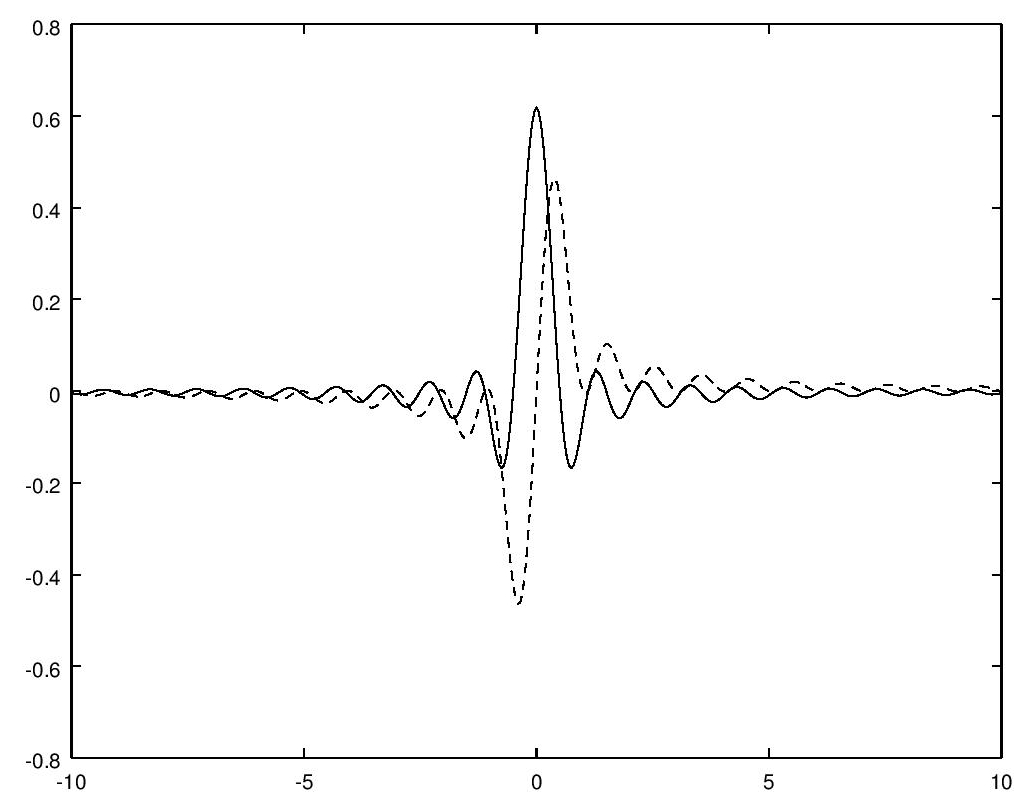}\hspace{0.5cm}
\includegraphics[scale=0.16]{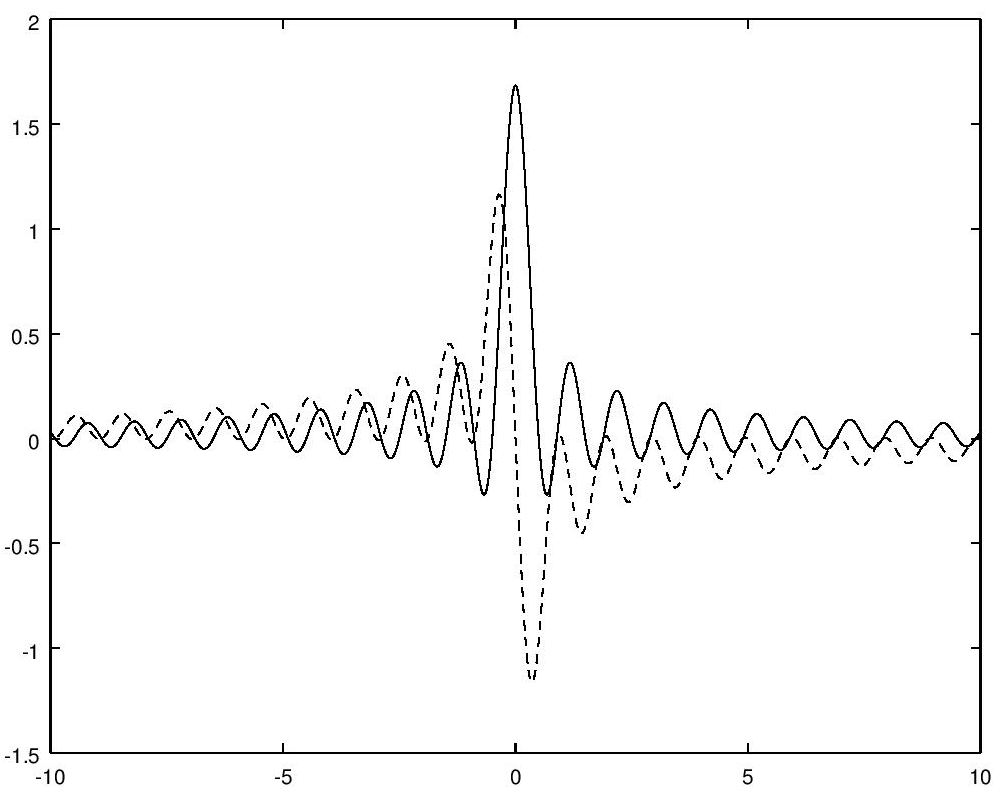}
\caption{Plot of the functions $g_1$ (left) and $\gamma_1$ (right) in Example \ref{ex}. The solid (resp. dashed) line shows the real (resp. imaginary) part of the functions $g_1$ and $\gamma_1$.\protect\footnotemark }
\label{plot_g_gamma}
\end{figure}
\end{example}

% and \textcolor{green}{ calculates the explicit values in Lemma %\ref{derivatives-nonzero} and Proposition }.
     
However, it turns out that already rather mild decay conditions on the time-frequency distribution of the windows exclude the possibility of reproducing pairs using two Gabor systems. Daubechies and Janssen \cite{dauja93} obtained a first result in this direction.
\begin{theorem}
Let  $\big(G(g,1,1),G(\gamma,1,1)\big)$ be a reproducing pair, then neither $g\in M_{2}^2$ nor $\gamma\in M^2_{2}$.
\end{theorem}     
In this paper we will show  a similar result where we replace the modulation space $M_{2}^2$ by $M_{1}^1$. Note that this is a new result as both spaces do not embed into each other.
     
     The following Proposition is a simple consequence of Janssen's  characterization of the modulation space $M^1$ via the Zak transform \cite{ja06}.
     \footnotetext{In the spirit of reproducible research we provide a Matlab/Octave script at \url{https://www.kfs.oeaw.ac.at/doc/RepPairGabor/rep_pair_gabor.m} which generates the plots of Figure \ref{plot_g_gamma}.}
\begin{proposition}\label{prop-modul-zak}
Let $s\in \N_0$ and $f\in M_{s}^1$, then $Z_af\in C^s(\R^2)$.
\end{proposition}
\proof Let $f\in M_{s}^1$ and $g_1,g_2\in M^1_{s}$ such that $Z_ag_1$ and
$Z_ag_2$ have no common zeros.  Adapting the argument of \cite[Theorem 4.1]{ja06} yields that $Z_af\cdot\overline{Z_ag_n},\ n=1,2,$ can be expressed by an  Fourier series 
$$
(Z_af\cdot\overline{Z_ag_n})(x,\omega)=\sum_{k,l\in\Z}c_n(k,l)e^{2\pi i kx/a+2\pi ial\omega}, \ n=1,2,
$$
with coefficients $\{c_n(k,l)\}_{k,l\in\Z}\in \ell^1_{v_s}(\Z^2)$.
Hence, if $|\beta|=\beta_1+\beta_2\leq s$, $\{k^{\beta_1}l^{\beta_2}c_n(k,l)\}_{k,l\in\Z}\in \ell^1(\Z^2)
$  which in turn implies that the Fourier series of the derivative
$$
D^\beta(Z_af\cdot\overline{Z_ag_n})(x,\omega)=\sum_{k,l\in\Z}(2\pi i)^{|\beta|}(k/a)^{\beta_1}(al)^{\beta_2}c_n(k,l)e^{2\pi i kx/a+2\pi ial\omega},
$$
converges absolutely for all $|\beta|\leq s$. Hence, $Z_af\cdot\overline{Z_ag_n} \in C^s(\R^2)$, for  $n=1,2$. We may choose $g_1,g_2$ to be Schwartz functions which guarantees that $Z_ag_n\in C^\infty(\R^2)$. As $Z_ag_1$ and $Z_ag_2$ have no common zeros, it follows that 
$
Z_af\in C^s(\R^2).
$ \pbox

\begin{lemma}\label{prop-no-schwartz}
Let $F\in L^2(\mathcal{Q}_a)$ be Lipschitz and assume that there exists $z^\ast\in \mathcal{Q}_a$ with $F(z^\ast)=0$, then $1/F\notin L^2(\mathcal{Q}_a)$. 
\end{lemma}
\proof Let $L$ be the Lipschitz constant of $F$ and $F(z^\ast)=0$, then 
$$
|F(z)|^2\leq L^2 \|z-z^\ast\|^2,\ \forall z\in B_\delta(z^\ast).
$$
Hence, $1/F\notin L^2([0,1]^2)$, as $1/|F(z)|^2\geq L^{-2}\|z-z^\ast\|^{-2}$ on $B_\delta(z^\ast)$.
\pbox

\begin{corollary}\label{new-balian-low}
Let 
$\big(G(g,a,1/a),G(\gamma,a,1/a)\big)$ be a reproducing pair,
then both $g\notin M_{1}^1$ and $\gamma\notin M_{1}^1$.
\end{corollary}
\proof Assume without loss of generality that $g\in M_{1}^1$, then 
$Z_ag\in C^1(\R^2)$ by Proposition \ref{prop-modul-zak}, $|Z_a\gamma|\geq m/|Z_ag|$ by \eqref{zak-cond-rep-pair} and there exists $z^\ast\in \mathcal{Q}_a$ such that $Z_a g(z^\ast)=0$. However, $Z_a\gamma\notin L^2(\mathcal{Q}_a)$ by Lemma \ref{prop-no-schwartz} and consequently $\gamma\notin L^2(\R)$, a contradiction.\pbox

\section{In quest of a reproducing partner for $\G$}\label{sec:quest}

For the rest of this paper, we focus on the study of  $\G$. As already mentioned in the introduction, $\G$ is a complete Bessel family but not a frame. By Corollary \ref{new-balian-low}, there is no dual Gabor system for $\G$. We will use Theorem \ref{theo-partner} to show that the expansion coefficients can be calculated via inner products. Finally, we show that any reproducing partner for $\G$ cannot have shift-invariant structure in time or frequency domain.

 \subsection{The range of $D_\G$}
 Condition $(i)$ in Theorem \ref{theo-partner} is satisfied. This is a consequence of \cite[Theorem 4.7]{ja81}, which states that 
 for every $f\in \S'(\R)$ there exists a sequence $\xi$ such that
$f=\sum_{n,m\in\Z}\xi[n,m]T_nM_m\varphi$ with  convergence
 in $\S'$-sense. Consequently, the series converges weakly for every $f\in L^2(\R)$ by density of 
 $\S(\R)$ in $L^2(\R)$, that is, $\ran D_\G=L^2(\R)$.\\

In order to verify condition $(ii)$, that is, that there exists $\{\xi_{k,l}\}_{(k,l)\in\Z^2}$ such that $D_\G\xi_{k,l}=e_{k,l}$ and $\sum_{k,l\in \Z}|\xi_{k,l}[n,m]|^2<\infty$, for all $(n,m)\in\Z^2$, we need some auxiliary results.
 \subsection{Solving $D_\G\xi=e_{k,l}$}

 \begin{lemma}\label{lemma1}
 Let $\gamma \in L^2(\R)$. The sequence $\xi_0$ is a (weak) solution of
  $D_\G \xi=\gamma$ 
  if, and only if, $\S_{k,l}\xi_0$ is a (weak) solution to 
  $D_\G \xi=\gamma_{k,l}$, where $S_{k,l}$ denotes the index-shift operator $$\S_{k,l}c[n,m]:=c[n-k,m-l].$$
 \end{lemma}
\proof Let $\xi_0$ be a (weak) solution to $D_\G \xi=\gamma$. It holds
\begin{align*}
\gamma_{k,l}&=T_{k}M_{l} \gamma=T_{k}M_{l} (T_\G \xi_0)
=T_{k}M_{l} \Big(\sum_{n,m\in\Z}\xi_0[n,m]T_{n}M_{m}\varphi\Big)\\ &=
\sum_{n,m\in\Z}\xi_0[n,m]T_{k+n}M_{l+m}\varphi =\sum_{n,m\in\Z}\xi_0[n-k,m-l]T_nM_{m}\varphi
\\ &=\sum_{n,m\in\Z}(\S_{k,l}\xi_0)[n,m]T_n M_{m}\varphi
=D_\G(\S_{k,l}\xi_0).
 \end{align*}
 The reversed implication follows with the same argument.
 \pbox

  Hence, in order to find all solutions of   $D_\G \xi=\gamma_{k,l}$, it remains to find one particular weak solution of $D_\G \xi=\gamma$ and characterize the kernel of $D_\G$ in $\V_\G(\Z^2)$.

\subsubsection{Characterizing the kernel of $D_\G$}

In this section, we will see that any weak solution of $D_\G\xi=0$ is given by a two dimensional polynomial  evaluated on $\Z^2$ times an oscillating sign factor. This result can already be found in \cite{ja81}.

As $\G$ is total in $L^2(\R)$ it follows that 
 $\xi\in \Ker D_\G$ if, and only if, $\ip{D_\G\xi}{\varphi_{n,m}}=0$, for every $(n,m)\in \Z^2$. It holds
\begin{align*}
\ip{D_\G\xi}{\varphi_{n,m}}&=\sum_{k,l\in\Z}\xi[k,l]\ip{T_k M_l \varphi}{T_n M_m \varphi}=
\sum_{k,l\in\Z}\xi[k,l]\ip{ T_{k-n}M_{l-m} \varphi}{\varphi}\\
&=\sum_{k,l\in\Z}\xi[k,l]\vartheta[n-k,m-l]=(\xi\ast \vartheta)[n,m],
\end{align*}
 where $$\vartheta[n,m]:=\ip{T_{-n}M_{-m}  \varphi}{\varphi}=(-1)^{ nm}e^{-\pi(n^2+m^2)/2},$$
 and the last equality holds by \cite[Lemma 1.5.2]{groe1}.

 \begin{comment}
 First we need the following result \cite[Lemma 1.5.2]{groe1}
 \begin{lemma}[Time-frequency shifts of Gaussians]\label{gauss-shift-lemma}
  For all $x,y,\omega,\eta\in\R$ it holds
  \begin{equation}\label{inner-prod-gauss}
   \ip{T_x M_\omega\varphi_\sigma}{T_yM_\eta\varphi_\sigma}=e^{\pi i(y-x)
   \cdot(\eta+\omega)}\varphi_{2\sigma}(x-y)\varphi_{2/\sigma}(\omega-\eta).
  \end{equation}
 \end{lemma}
For integer time-frequency shifts and $\sigma=1$, equation \eqref{inner-prod-gauss} yields
$$
   \ip{T_n M_m\varphi}{T_kM_l\varphi}=e^{\pi i(n-k)\cdot(m-l)}\varphi_{2}(n-k)
   \varphi_{2}(m-l),
$$

 Let us consider the sequence $\vartheta(n,m):=\ip{T_{-n} M_{-m}\varphi}{\varphi}$. Since $\G(\gamma,1,1)$ is
 an orthonormal basis, it follows
\begin{align*}
 \vartheta(n,m)&=
 \sum_{k,l\in\Z}\Big\langle T_{-n}M_{-m}\varphi,\ip{\varphi}{T_kM_l\gamma} T_kM_l\gamma\Big\rangle\\
 &=
 \sum_{k,l\in\Z}\ip{T_kM_l\gamma}{\varphi}\ip{ T_{-n}M_{-m}\varphi}{ T_kM_l\gamma}\\
  &=
 \sum_{k,l\in\Z}\ip{\gamma}{T_{-k}M_{-l}\varphi}\ip{T_{-n-k} M_{-m-l}\varphi}{\gamma}\\
   &=
 \sum_{k,l\in\Z}\ip{\gamma}{M_{l}T_{k}\varphi}\ip{ T_{-(n-k)}M_{-(m-l)}\varphi}{\gamma}\\
 &= \sum_{k,l\in\Z}\overline{a(-k,-l)} a(n-k,m-l)
 =(a^\ast \ast a)(n,m),
\end{align*}
  where $a^\ast(n,m):=\overline{a(-n,-m)}$ denotes the involution operator.
  \end{comment}

 Let  $\xi\ast \vartheta=0$, then $\F_d(\xi)\cdot\F_d(\vartheta)=0$ holds at least in the sense of periodic distributions. Hence, 
 we intend to characterize those periodic distributions $\Lambda$ such that $\Lambda\cdot \Theta=0$, where   $\Theta:=\F_d(\vartheta)$.
 The function $\Theta$ can be expressed analytically as
  \begin{equation}\label{theta-fourier}
  \Theta(\omega)= \sum_{n,m\in\Z}
 (-1)^{ nm}e^{-\pi(n^2+m^2)/2}e^{-2\pi i(n\omega_1+m\omega_2)}
\end{equation}
 $$
  =\theta_3(\pi\omega_1,e^{-\pi/2})\cdot\theta_3(2\pi\omega_2,e^{-2\pi})+
  \theta_4(\pi\omega_1,e^{-\pi/2})\cdot\theta_2(2\pi\omega_2,e^{-2\pi}),$$
where $\omega\in[0,1)^2$ and $\theta_k$ denotes the $k$-th Jacobi theta function, see for example \cite{whiwat96}. 
This function has been studied in 
\cite[Theorem 3.5]{ja81}. We will state the results in the following Lemma.

\begin{comment}Throughout the rest of this paper we will make use of multi-indices to denote partial derivatives. Let $\alpha,\beta\in \N_0^d$ with $\alpha_n\geq\beta_n$ for all $n=1,...,d$, then  $D^\alpha F:=\dfrac{\partial^{|\alpha|}}{\partial_{x_1}^{\alpha_1}\ldots \partial_{x_d}^{\alpha_d}}$, where $|\alpha|:=\alpha_1+...+\alpha_d$. 
\end{comment}
     \begin{lemma}\label{derivatives-nonzero}
     Set $\omega_0:=(1/2,1/2)$. The function $\Theta\in C^\infty([0,1]^2)$ has the following properties:
     \begin{enumerate}[(i)]
      \item $\Theta(\omega)\geq0,\ \forall\ \omega\in[0,1)^2$, 
      \item
  $\Theta(\omega)=0$ if, and only if, $\omega=\omega_0$,
  \item  $D^{(2,0)}\Theta(\omega_0)=D^{(0,2)}\Theta(\omega_0)> 0$.
  \begin{comment}
    \item whenever $k$ or $l$ is odd, then $
   D^{(k,l)}\Theta(\omega_0)=0,$
  \item for every $k\in2\N_0$, it holds
  $$  D^{(k,2)}\Theta(\omega_0)=D^{(2,k)}\Theta(\omega_0)\neq 0,
  $$
  and  
  $$
 D^{(k,4)}\Theta(\omega_0)=D^{(4,k)}\Theta(\omega_0)\neq 0.
  $$
  \end{comment}
     \end{enumerate}
 \end{lemma}
 \begin{comment}
 \proof The  statements $(i)$ - $(iv)$ were shown in \cite[Theorem 3.5]{ja81}. It remains to prove $(v)$.
 It follows directly from \eqref{theta-fourier} that the partial derivatives of $\Theta$ at $\omega_0$ can be written as 
   $$
D^{(k,l)}\Theta(\omega_0)=
  (-2\pi i)^{k+l}\sum_{n,m\in\Z} n^k m^l e^{-\pi(n^2+m^2)/2}(-1)^{n+m+nm}.
  $$
  If $m$ is even, then 
  $$(-1)^{n+m+nm}=\left\{\begin{array}{ll} 1, &\mbox{ for n  even}\\-1, &\mbox{ for n  odd.}
                        \end{array}\right.
$$ If $m$ is odd, then $(-1)^{n+m+nm}=-1$ 
 for every $n\in\Z$. Hence,
    $$
    (-2\pi i)^{-(k+l)} D^{(k,l)}\Theta(\omega_0)
  $$
  $$
  =
  \sum_{n\in\Z}(-1)^n n^k e^{-\pi n^2/2}  \sum_{m\in2\Z}  m^l  e^{-\pi m^2/2}
  -  \sum_{n\in\Z} n^k e^{-\pi n^2/2}  \sum_{m\in2\Z+1} m^l e^{-\pi m^2/2},
  $$
which implies that the partial derivative of $\Theta$ at $\omega_0$ is zero if, and only if,
   $$
 Q(k):=\dfrac{\sum_{n\in\Z} (-1)^n n^k  e^{-\pi n^2/2}}{\sum_{n\in\Z} n^k e^{-\pi n^2/2}}
=\dfrac{ \sum_{m\in2\Z+1} m^l e^{-\pi m^2/2} }{ \sum_{m\in2\Z} m^l e^{-\pi m^2/2}  }=:P(l).
  $$
  The value $P(2)$ (resp. $P(4)$) can be calculated numerically and equals approximately $27.83$ (resp. 6.96). On 
  the other hand, $|Q(k)|\leq 1,$  for every $k\in 2\N_0$ which concludes the proof.\pbox
  \end{comment}
  
  A particular consequence of Lemma \ref{derivatives-nonzero} is that any periodic distribution $\Lambda$ satisfying 
  $\Lambda\cdot\Theta=0$ is supported on $\{\omega_0\}$.
Thus, by \cite[Theorem 6.25]{ru91} there exists $N\in\N_0$ and coefficients $c_\alpha\in \C$, $\alpha\in\N^2_0$, such that
   \begin{equation}\label{Lambda}\Lambda=\sum_{|\alpha|\leq N} c_\alpha D^\alpha\delta_{\omega_0},\end{equation}

   Applying the inverse Fourier transform to $\Lambda$  immediately shows the 
   following result.
  \begin{corollary}\label{sequence-shape-corr}
   Every sequence  $p\in\Ker D_\G$ can be written as 
   $$
   p[n,m]=(-1)^{n+m}\sum_{|\alpha|\leq N}c_{\alpha}\cdot n^{\alpha_1} m^{\alpha_2}.
   $$
   \begin{comment}
 where the coefficients $c_{k,l}$ satisfy the following conditions
\begin{align*}
    &(a)\ c_{2,0}+c_{0,2}=0 &(c)\ 3c_{0,3}+c_{2,1}=0 &\hspace{0.7cm}(e)\ c_{4,0}= c_{0,4}=c_{2,2}=0.\\
   &(b)\ c_{3,1}+c_{1,3}=0 &(d)\ 3c_{3,0}+c_{1,2}=0 &\mbox{}
\end{align*}
      In particular,
      $$
   p[n,0]=(-1)^n\sum_{k\leq 3}c_{k,0}\cdot n^k .
      $$
      \end{comment}
  \end{corollary}

      \subsubsection{Calculation of $\xi_0$}
      We choose the orthonormal basis $\{e_{n,m}\}_{n,m\in\Z}$ in Theorem \ref{theo-partner} to be the Gabor system 
  $G(\gamma,1,1)$ with $\gamma:=\chi_{[-1/2,1/2]}$.
By Lemma \ref{lemma1} it remains to find a weak solution of
   $D_\G \xi=\gamma$.\\
      
      A first attempt in finding the expansion coefficients for $\G$ can be found in \cite{ba80}. In this paper, Bastiaans constructed a window that yields the expansion coefficients
for a dense subspace of $L^2(\R)$, see also \cite[Section 4.4]{ja82}.
   Bastiaans' dual window  \cite{ba80} is analytically given by 
 \begin{equation}\label{bastiaans_dual}
    \psi(t)=C_\psi\ e^{\pi t^2}\sum_{n>|t|-1/2}(-1)^{n}e^{-\pi(n+1/2)^2},
   \end{equation}
    for some constant $C_\psi>0$. The function $\psi$ is defined as 
    $\psi:=Z^{-1}(1/\overline{Z(\varphi)})$ and has the property that $\psi$ is bounded but not contained in
    $L^p(\R)$ for any $1\leq p<\infty$, that is, $\psi\in L^\infty(\R)\backslash L^p(\R)$.

    Janssen showed in \cite[Section 4.4]{ja82} that, for every $f\in L^1(\R)$ satisfying $|\ip{f}{T_k M_l\psi}|\rightarrow 0$,
    for every fixed $l\in\Z$ and $|k|\rightarrow \infty$, the following series converges in the sense of tempered distributions
    \begin{equation}\label{bastiaans-exp}
    f=\sum_{k,l\in\Z}\ip{f}{T_k M_l\psi}T_kM_l\varphi.
    \end{equation}
We will prove that this condition implies weak convergence for a dense subspace of $L^2(\R)$.
    \begin{lemma}\label{dense_subspace_coefficients}
     For every function $f\in \mathcal{M}$, where
     $$
     \mathcal{M}:=\Big\{f\in L^2(\R):\ f=\sum_{k,l\in\Z}c[k,l]T_kM_l\gamma,\ \|c\|_1<\infty\Big\},
     $$
    equality \eqref{bastiaans-exp} holds weakly in $L^2(\R)$.
  In particular, $\xi_0[k,l]:=\ip{\gamma}{T_k M_l\psi}$ is a weak solution of $T_\G\xi=\gamma$. 
    \end{lemma}
\proof 
 Recall that  $G(\gamma,1,1)$ is an orthonormal basis. Let $f\in\mathcal{M}$, then $f\in L^1(\R)\cap L^2(\R)$ as
 $$
 \|f\|_1\leq\sum_{k,l\in\Z}|c[k,l]|\|T_kM_l\gamma\|_1=\|c\|_1\ \mbox{ and }\ 
 \|f\|_2^2=\sum_{k,l\in\Z}|c[k,l]|^2\leq \|c\|_1^2.
 $$
 
 Let $k\neq 0$, then
  \begin{align*}
    |\ip{\gamma}{T_k M_l\psi}|&\leq \int_{-1/2}^{1/2}|\psi(t-k)|dt\leq C \int_{-1/2}^{1/2}e^{\pi (t-k)^2}dt\sum_{n\ge|k|}
    e^{-\pi(n+1/2)^2}
  \\
    &=C \int_{-1/2}^{1/2}e^{\pi (t-|k|)^2-\pi(|k|+1/2)^2}dt\sum_{n\geq|k|}e^{-\pi(n+1/2)^2+\pi(|k|+1/2)^2}\\
    &\leq C\int_{-1/2}^{1/2}e^{\pi (t^2-1/4-|k|(2t+1))}dt\leq C\int_{-1/2}^{1/2}e^{-\pi |k|(2t+1)}dt\\
    &=\frac{C}{2\pi |k|}(1-e^{-2\pi |k|})\leq \frac{C}{2\pi(1+ |k|)},
    \end{align*}
    where we have used that $e^{\pi t^2}$ is even and that  $t^2-1/4<0$ on $(-1/2,1/2)$.
     Observe that the constant $C$ is independent of 
    $k$ as 
    $$\sum_{n\geq|k|}e^{-\pi(n+1/2)^2+\pi(|k|+1/2)^2}=\sum_{n\geq0}e^{-\pi(n+|k|+1/2)^2+\pi(|k|+1/2)^2}$$
    $$
    =\sum_{n\geq0}e^{-\pi(n^2+2n|k|+n)}\leq\sum_{n\geq0}e^{-\pi(n^2+n)}=C<\infty.
    $$
    Let now $f=\sum_{n,m}c[n,m]T_nM_m \gamma$ with $c\in\ell^1(\Z^2)$.
    By Young's inequality for convolutions and previous calculations we get
    \begin{align*}
     \sum_{k\in \Z}|\ip{f}{T_kM_l\psi}|^2&=\sum_{k\in \Z}\Big|\sum_{n,m\in\Z}c[n,m]\ip{\gamma}{T_{k-n}M_{l-m}\psi}\Big|^2\\
     &\leq C\sum_{k\in \Z}\Big(\sum_{n,m\in\Z}|c[n,m]|(1+|k-n|)^{-1}\Big)^2\\
     &= C\sum_{k\in \Z}\Big(\sum_{n\in\Z}\Big(\sum_{m\in\Z} c[n,m]\Big)(1+|k-n|)^{-1}\Big)^2\\
     &\leq C\|c\|_1^2\sum_{k\in\Z} (1+|k|)^{-2}<\infty.
    \end{align*}
This estimate finally implies that for fixed $l\in\Z$, $|\ip{f}{T_kM_l\psi}|\rightarrow 0$, as $|k|\rightarrow \infty$.
\pbox
      
Let us investigate the behavior of $\xi_0$ in more detail. Thereby, we will use the notation $g(t):=C_\psi e^{\pi t^2}$,
$$G_k:=\sum_{n\geq0}(-1)^{n}
    e^{-\pi(n^2+2|k|n+n+1/4)},$$
     and $\mu_k(t):=(-1)^ke^{-\pi|k|}G_k e^{-2\pi kt}$.
It holds, 
\begin{align}
\xi_0[k,l]&=\ip{\gamma}{T_k M_l\psi}=C_\psi\int_{-1/2}^{1/2}e^{\pi(t-k)^2}e^{-2\pi i lt}dt\sum_{n\geq|k|}(-1)^n
e^{-\pi(n+1/2)^2} \nonumber \\
&=C_\psi\int_{-1/2}^{1/2}e^{\pi t^2}e^{-2\pi kt}e^{-2\pi i lt}dt\sum_{n\geq|k|}(-1)^n
 e^{-\pi(n+1/2)^2}e^{\pi k^2}\nonumber \\
&=C_\psi\int_{-1/2}^{1/2}e^{\pi t^2}e^{-2\pi kt}e^{-2\pi i lt}dt\sum_{n\geq0}(-1)^{n+k}
    e^{-\pi(n^2+2|k|n+|k|+n+1/4)}\nonumber \\
 &=C_\psi(-1)^ke^{-\pi|k|}G_k\int_{-1/2}^{1/2}e^{\pi t^2}e^{-2\pi kt}e^{-2\pi i lt}dt\nonumber \\
 &=\int_{-1/2}^{1/2}g(t)\mu_k(t)e^{-2\pi i lt}dt=\mathcal{F}(g\cdot \mu_k)[l]. \label{xi_0=f(gmu)}
\end{align}

For $k=0$, the Fourier transform of $\mu_0$ is given by 
\begin{equation}\label{Fmuk=n}
\F(\mu_0)[l]=G_0\delta_0[l].
\end{equation}
Let $k\neq 0$, then
\begin{align}
\F(\mu_k)[l]&=(-1)^ke^{-\pi|k|}G_k\int_{-1/2}^{1/2}e^{-2\pi (k+ i l)t}dt\nonumber \\
&=\frac{(-1)^{k+1}e^{-\pi|k|}G_k}{2\pi (k+il)}(e^{-\pi (k+ i l)}-e^{\pi (k+ i l)})\nonumber \\
&=\frac{(-1)^{l+k+1}}{2\pi (k+il)}(e^{-\pi(|k|+k)}-e^{-\pi(|k|- k)})G_k\nonumber \\
&=\frac{(-1)^{l+k}}{2\pi (k+il)}\text{sgn}(k)(1-e^{-2\pi|k|})G_k
\nonumber\\&=\frac{(-1)^{l+k}}{2\pi (k+il)}H_k,\label{F(mu_k)}
\end{align}
where $H_k:=\text{sgn}(k)(1-e^{-2\pi|k|})G_k$.
Observe that $H_k\rightarrow\pm e^{-\pi/4}$ as $k\rightarrow\pm\infty$.

  \subsection{Existence of reproducing partners}\label{sec:existence}
  
  \begin{comment}
The following definition will provide a necessary condition for the convergence of a series in $\Z^2$.
   \begin{definition}
We say that a subset $I\subset \Z^2$ possesses property $(D)$ if
\begin{equation}\label{property-p}
\dfrac{|I\cap  B_r^\infty|}{|B_r^\infty|}\rightarrow 1,\ \mbox{as } r\rightarrow \infty,
\end{equation}
 where  $B_r^\infty:=\{(k,l)\in\Z^2:\ \max(|k|,|l|)\leq r\}$.
 \end{definition}
 Note that, if a finite collection of index sets
 $I_1,I_2,...,I_N\subset\Z^2$ satisfies $(D)$,
then  $\bigcap_{n=1}^NI_n$ also satisfies $(D)$. To see this, let $I_1$ and  $I_2$ possess property $(D)$, then $$\dfrac{|B_r^\infty\cap (I_1\cap I_2)^c|}{|B_r^\infty|} =\dfrac{|(B_r^\infty\cap I_1^c)\cup (B_r^\infty\cap I_2^c)|}{|B_r^\infty|}\leq \dfrac{|B_r^\infty\cap I_1 ^c|}{|B_r^\infty|}+\dfrac{|B_r^\infty\cap I_2^c|}{|B_r^\infty|}\rightarrow 0,$$ which implies that $(D)$ holds for $I_1\cap I_2$.
 
 \begin{remark}\label{convergence-remark}
Let $\{a_{k,l}\}_{k,l\in\Z}\subset \R_+$, if the series 
$$
\sum_{(k,l)\in\Z^2}\frac{a_{k,l}}{1+k^2+l^2}
$$
converges, then it follows that the index set 
$$I(\varepsilon):=\{(k,l)\in\Z^2:\ a_{k,l}<\varepsilon\}$$ satisfies $(D)$, for every $\varepsilon>0$.
\end{remark}    
We are now able to prove our main result. We will thereby use an argument based on Remark \ref{convergence-remark}.
\end{comment}

We are now ready to state and prove our main result.
  \begin{theorem}\label{main-repr}
   There exists a system $\Psi$ making $(\Psi,\G)$ a reproducing pair. In other words, there exists a dual system for the
   complete Bessel sequence $\G$.
  \end{theorem}
  
\proof  Set 
\begin{equation}\label{xikl=Sxi0+pkl}\xi_{k,l}=\mathcal{S}_{k,l}\xi_0+p_{k,l},
\end{equation}
  where $p_{k,l}\in \Ker D_\G$, then $D_\G \xi_{k,l}=\gamma_{k,l}$. Let us assume 
  that there exist a reproducing partner for $\G$. By  Theorem \ref{theo-partner}, we can choose $\{\xi_{k,l}\}_{k,l\in\Z}$ in such a way that
\begin{equation}\label{whole-sum-conv-assumpt}
  \sum_{k,l\in\Z}\big|\xi_{k,l}[n,m]\big|^2<\infty,\ \forall\ n,m\in\Z.
  \end{equation}
Recall that,  by Corollary \ref{sequence-shape-corr}, $p_{k,l}$ is given by $$p_{k,l}[n,m]=(-1)^{n+m}\sum_{|\alpha|\leq N}c_{\alpha}[k,l]\cdot n^{\alpha_1}m^{\alpha_2}.$$
In the following we will  choose $N=0$ for every $(k,l)\in\Z^2$, that is, $p_{k,l}[n,m]=(-1)^{n+m}c[k,l]$.
As we assume \eqref{whole-sum-conv-assumpt}, we can apply Parseval's formula to the summation with respect to
$l\in \Z$. Using  \eqref{xi_0=f(gmu)} thus yields
\begin{align}\label{estimates123}
  \sum_{l\in\Z}\big|\xi_{k,l}&[n,m]\big|^2=
  \sum_{l\in\Z}\big|\S_{k,l}\xi_0[n,m]+p_{k,l}[n,m]\big|^2\nonumber \\
  &=  \sum_{l\in\Z}\big|\xi_0[n-k,m-l]+p_{k,l}[n,m]\big|^2\nonumber\\  
  &=\sum_{l\in\Z}\Big|
  \F(g\cdot \mu_{n-k})[m-l]+(-1)^{n+m}c[k,l]\Big|^2 \nonumber
  \\ &=\int_{-1/2}^{1/2}\Big|M_{m}(g\cdot\mu_{n-k})(-\omega) +(-1)^{n+m}\F^{-1}(c[k,\cdot])(\omega)\Big|^{2}d\omega \nonumber
\\
&\leq |g(1/2)|^2\int_{-1/2}^{1/2}\Big|M_{m}\mu_{n-k}(-\omega) +(-1)^{n+m}
(\F^{-1}(c[k,\cdot])/g)(\omega)\Big|^{2}d\omega \nonumber
\\
&=C\sum_{l\in\Z}\Big|\F(\mu_{n-k})[m-l]+(-1)^{n+m}\F\big(\F^{-1}(c[k,\cdot])/g\big)[l]\Big|^{2}=:(\ast).
\end{align}
If $k\neq n$ it follows by \eqref{F(mu_k)} that
\begin{align}
  (\ast)&= C\sum_{l\in\Z}\Big|\frac{(-1)^{n+k+l+m}H_{n-k}}{2\pi(n-k+i(m-l))}+(-1)^{n+m}\F\big(\F^{-1}(c[k,\cdot])/g\big)[l]\Big|^{2} \nonumber
\\
&= \frac{C}{4\pi^2}\sum_{l\in\Z}\frac{\big|H_{n-k}-
\beta[k,l]\cdot\big(1-\frac{n+im}{k+il}\big)\big|^{2}}{(n-k)^2+(m-l)^2} \label{sum-lower-bound-xi},
\end{align}
where
\begin{equation*}
\beta[k,l]:=2\pi(-1)^{k+l}\F\big(\F^{-1}(c[k,\cdot])/g\big)[l]\cdot(k+il).
\end{equation*}
If $k=n$, then by \eqref{Fmuk=n} and \eqref{estimates123} we get the estimate
\begin{equation}\label{sum-if-n=k}
  \sum_{l\in\Z}\big|\xi_{n,l}[n,m]\big|^2\leq C \sum_{l\in\Z}|G_0\delta_0[m-l]-2\pi^{-1}(-1)^{m+l}\beta[n,l]/(n+il)|^2.
\end{equation}
Let us choose $c_{k,l}$ such that $c_{0,l}:=0$ for every $l\in\Z$ and  
$$
c[k,l]:=(2\pi)^{-1}\text{sgn}(k)e^{-\pi/4}\cdot\F(\F^{-1}(h_k)\cdot g)[l],\ \mbox{if } k\neq 0,
$$
where $h_k[l]:=(-1)^{k+l}(k+il)^{-1}$.
 Then  $\beta[k,l]=\text{sgn}(k)e^{-\pi/4}$ and therefore the right hand side of \eqref{sum-if-n=k} converges for every $(n,m)\in \Z^2$. Moreover, we have
$$
  \sum_{\substack{k,l\in\Z\\ k\neq n}}\big|\xi_{k,l}[n,m]\big|^2
  $$
  $$
   \leq C  \sum_{\substack{k,l\in\Z\\ k\neq n}}\frac{\big|H_{k-n}-\text{sgn}(k)e^{-\pi/4}\big|^2}{(n-k)^2+(m-l)^2}+\frac{e^{-\pi/4}\cdot(1-\delta_0[k])\cdot(n^2+m^2)}{(k^2+l^2)\cdot\big((n-k)^2+(m-l)^2\big)}.
$$
It is easy to see that the second series is finite and the first one can be estimated as follows
$$
\sum_{\substack{k,l\in\Z\\ k\neq n}}\frac{\big|H_{k-n}-\text{sgn}(k)e^{-\pi/4}\big|^2}{(n-k)^2+(m-l)^2}
\leq \sum_{\substack{k\in\Z\\ k\neq n}}\big|H_{k-n}-\text{sgn}(k)e^{-\pi/4}\big|^2\sum_{l\in\Z}\frac{1}{1+l^2},
$$
which is finite for every $(n,m)\in\Z^2$. 
All in all, we have shown that \eqref{whole-sum-conv-assumpt} holds for every $(n,m)\in\Z^2$. This concludes the proof.
\pbox

\subsection{Non-existence of shift-invariant reproducing partners}
We conclude this paper by showing that any reproducing partner for $\G$ is necessarily unstructured in the sense that it cannot be written as a shift-invariant system neither in time nor in frequency domain.
\begin{proposition}
There exists no shift-invariant dual system for $\G$, that is, any reproducing partner $\Psi$ for $\G$ cannot be written as $\psi_{n,m}=T_n\psi_m$ or $\psi_{n,m}=M_m\psi_n$.
\end{proposition}
  \proof Let us write $\psi_{n,m}$ as a Gabor expansion with respect to the orthonormal basis $G(\gamma,1,1)$ and assume that $\psi_{n,m}=T_n\psi_m$, then 
\begin{align*}
  \psi_{0,m}&=T_0\psi_m=T_{-n}\psi_{n,m}=\sum_{k,l\in\Z}\overline{\xi_{k,l}[n,m]}T_{k-n}M_l\gamma
\\
  &=\sum_{k,l\in\Z}\overline{\xi_{k+n,l}[n,m]}T_{k}M_l\gamma,
\end{align*}
  which implies that 
  $\xi_{k,l}[0,m]=\xi_{k+n,l}[n,m]$, for all $(n,m,k,l)\in\Z^4$.
Since, by \eqref{xikl=Sxi0+pkl},
  $$
  \xi_{k+n,l}[n,m]=\S_{k+n,l}\xi_0[n,m]+p_{k+n,l}[n,m]$$
$$  
  =\xi_0[-k,m-l]+p_{k+n,l}[n,m]=\S_{k,l}\xi_0[0,m]+p_{k+n,l}[n,m],
  $$
  it thus follows that
   \begin{equation}\label{eq-for-pol}
   p_{k,l}[0,m]=p_{k+n,l}[n,m],\ \forall\ n,m,k,l\in\Z.
   \end{equation}
  In the following, we will denote $c_0:=c_{0,0}$ and let $\alpha,\beta\in\N_0^2$.  By Corollary \ref{sequence-shape-corr} and \eqref{eq-for-pol}, one has
  $$
  (-1)^{m}\sum_{\substack{|\beta|\leq N\\ \beta_1=0}} c_{\beta}[k,l] \cdot m^{\beta_2}=(-1)^{n+m}\sum_{|\alpha|\leq N} c_{\alpha}[k+n,l]\cdot n^{\alpha_1} m^{\alpha_2}.
  $$
Setting $k=s-n$, $m=0$  then yields 
  $$
c_{0}[s-n,l]=(-1)^n\sum_{\substack{|\alpha|\leq N\\ \alpha_2=0}} c_{\alpha}[s,l]\cdot n^{\alpha_1}.
  $$
If there exists $(s,l)\in\Z^2$, such that 
\begin{equation}\label{c_alpha-comparison}\sum_{\substack{|\alpha|\leq N\\ \alpha_2=0}} c_{\alpha}[s,l]\cdot n^{\alpha_1}\neq 0,\ \mbox{ for some }n\in\Z,
\end{equation}
then either $|c_{0}[s-n,l]|\rightarrow \infty$, as $|n|\rightarrow\infty$, or $c_{0}[s-n,l]=(-1)^{n}c_{0}[s,l]\neq 0$, for all $n\in\Z$, which implies $|c_{0}[\cdot,l]|\equiv C>0$. As $|\xi_0[k,l]|\rightarrow 0$, for $|(k,l)|\rightarrow \infty$, we obtain for both cases that
\begin{equation}\label{sum-for-00-infty}\sum_{k,l\in\Z}|\xi_{k,l}[0,0]|^2=\infty.
\end{equation}
If there exist no $(s,l)\in\Z^2$, such that \eqref{c_alpha-comparison} holds, then $p_{k,l}[0,0]\equiv 0$ and consequently  $\xi_{k,l}[0,0]=\xi_0[-k,-l]$. Hence, \eqref{sum-for-00-infty} holds and $\psi_{0,0}$ is not well defined.
\\ An analogous argument can be used to show that $\psi_{n,m}$ cannot be given by $M_m\psi_n$. \pbox
  
  \section{Conclusion}
    It appears that the obstructions to Gabor frames as portrayed for example in \cite{groe14} are preserved if one leaves the setup of frame theory and considers reproducing pairs consisting of two Gabor families instead. Already minor decay conditions on the window functions exclude the possibility to form reproducing pairs. \\
   We have seen that the concept of reproducing pairs provides new insights on complete vector systems. In particular, the characterization given in Theorem \ref{theo-partner} has proven to be a useful tool for different vector families, see also \cite{ansptr15}.\\
  Here, we used Theorem \ref{theo-partner} to show the existence of dual systems for the system of integer time-frequency shifts of the Gaussian. The crucial point in our argument is to estimate the behaviour of $\xi_0$ and choose appropriate elements of the kernel of $D_\G$. We believe that the same recipe will work for other windows, like the Hermite functions.

    \section*{Acknowledgement}
This work was funded by the Austrian Science Fund (FWF) START-project FLAME ('Frames and
Linear Operators for Acoustical Modeling and Parameter Estimation'; Y 551-N13).
The authors would like to thank Tomasz Hrycak and Hartmut F{\"u}hr for valuable discussions.

\bibliographystyle{plain}
\bibliography{paperbib}

\begin{thebibliography}{10}

\bibitem{algro01}
A.~Aldroubi and K.~Gr{\"o}chenig.
\newblock Non-uniform sampling and reconstruction in shift-invariant spaces.
\newblock {\em SIAM Review}, 43:585--620, 2001.

\bibitem{alanga00}
S.~T. Ali, J-P. Antoine, and J-P. Gazeau.
\newblock {\em Coherent {S}tates, {W}avelets and their {G}eneralizations}.
\newblock Graduate texts in contemporary physics. Springer, 2nd edition, 2014.

\bibitem{jpaxxl09}
J-P. Antoine and P.~Balazs.
\newblock Frames and semi-frames.
\newblock {\em J. Phys. A: Math. Theor.}, 44(20):205201, 2011.

\bibitem{jpaxxl12}
J-P. Antoine and P.~Balazs.
\newblock Frames, semi-frames and {H}ilbert scales.
\newblock {\em Numer. Func. Anal. Opt.}, 33:736--769, 2012.

\bibitem{ansptr15}
J-P. Antoine, M.~Speckbacher, and C.~Trapani.
\newblock Reproducing pairs of measurable functions.
\newblock {\em submitted to: Acta Appl. Math.}, 2015.
\newblock arXiv:1505.04187v2.

\bibitem{bagrza75}
H.~Bacry, A.~Grossmann, and J.~Zak.
\newblock Proof of completeness of lattice states in the kq representation.
\newblock {\em Phys. Rev. B}, 12(4), 1975.

\bibitem{babugikl71}
V.~Bargmann, P.~Butera, L.~Girardello, and J.~R. Klauder.
\newblock On the completeness of the coherent states.
\newblock {\em Rep. Math. Phys.}, 2(4):221--228, 1971.

\bibitem{ba80}
M.~J. Bastiaans.
\newblock Gabor's expansion of a signal into {G}aussian elementary signals.
\newblock {\em Proc. IEEE}, 68(4):538--539, 1980.

\bibitem{christ1}
O.~Christensen.
\newblock {\em An {I}ntroduction to {F}rames and {R}iesz {B}ases.}
\newblock {A}pplied and {N}umerical {H}armonic {A}nalysis. {B}irkh{\"a}user,
  2003.

\bibitem{dauja93}
I.~Daubechies and A.~J. E.~M. Janssen.
\newblock Two theorems on lattice expansions.
\newblock {\em IEEE Trans. Inform. Theory}, 39(1):3--6, 1993.

\bibitem{duscha52}
R.~J. Duffin and A.~C. Schaeffer.
\newblock A class of nonharmonic {F}ourier series.
\newblock {\em Trans. Amer. Math. Soc.}, 72:341--366, 1952.

\bibitem{hgfei06}
H.~G. Feichtinger.
\newblock Modulation spaces: Looking back and ahead.
\newblock {\em Samp. Theory Signal Image Process.}, 5:109--140, 2006.

\bibitem{fezi98}
H.~G. Feichtinger and G.~Zimmermann.
\newblock A {B}anach space of test functions for {G}abor analysis.
\newblock In H.G. Feichtinger and T.~Strohmer, editors, {\em Gabor analysis and
  {A}lgorithms - {T}heory and {A}pplications}. Birkh{\"a}user Boston, 1998.

\bibitem{ga46}
D.~Gabor.
\newblock Theory of communication. {P}art 1: {T}he analysis of information.
\newblock {\em J. Inst. Electr. Eng. 3}, 93:429--441, November 1946.

\bibitem{groe1}
K.~Gr{\"o}chenig.
\newblock {\em {F}oundations of {T}ime-{F}requency {A}nalysis}.
\newblock {A}ppl. {N}umer. {H}armon. {A}nal. {B}irkh{\"a}user {B}oston, 2001.

\bibitem{groe14}
K.~Gr{\"o}chenig.
\newblock The mystery of {G}abor frames.
\newblock {\em J. Fourier Anal. Appl.}, 20(4):865--895, 2014.

\bibitem{hepo06}
C.~Heil and M.~Powell.
\newblock Gabor {S}chauder bases and the {B}alian-{L}ow theorem.
\newblock {\em J. Math. Phys.}, 74, 2006.

\bibitem{ja81}
A.~J. E.~M. Janssen.
\newblock Gabor representation of generalized functions.
\newblock {\em J. Math. Anal. Appl.}, 83:377--394, 1981.

\bibitem{ja81a}
A.~J. E.~M. Janssen.
\newblock Weighted {W}igner distributions vanishing on lattices.
\newblock {\em J. Math. Anal. Appl.}, 80:156--167, 1981.

\bibitem{ja82}
A.~J. E.~M. Janssen.
\newblock Bargmann transform, {Z}ak transform, and coherent states.
\newblock {\em J. Math. Phys.}, 23(5):720--731, 1982.

\bibitem{ja06}
A.~J. E.~M. Janssen.
\newblock Zak transform characterization of ${S}_0$.
\newblock {\em Samp. Theory Signal Image Process.}, 5(2):141--162, 2006.

\bibitem{pe71}
A.~M. Perelomov.
\newblock On the completeness of a system of coherent states.
\newblock {\em Theoret. Math. Phys.}, 6(2):156--164, 1971.

\bibitem{ru91}
W.~Rudin.
\newblock {\em Functional analysis}.
\newblock International Series in Pure and Applied Mathematics. McGraw-Hill
  Inc., New York, second edition, 1991.

\bibitem{spexxl14}
M.~Speckbacher and P.~Balazs.
\newblock Reproducing pairs and the continuous nonstationary {G}abor transform
  on {LCA} groups.
\newblock {\em J. Phys. A: Math. Theor.}, 48(395201), 2015.

\bibitem{neu55}
J.~von Neumann.
\newblock {\em Mathematische {G}rundlagen der {Q}uantenmechanik}.
\newblock Springer, Berlin, 1932.
\newblock English translation: "{M}athematical {F}oundations of {Q}uantum
  {M}echanics," Princeton Univ. Press, 1955.

\bibitem{whiwat96}
E.~T. Whittaker and G.~N. Watson.
\newblock {\em A {C}ourse of {M}odern {A}nalysis}.
\newblock Cambridge University Press, 4th edition, 1996.

\end{thebibliography}
 
\end{document}